\documentclass[a4paper,12pt]{article}

\usepackage{amsmath,amsfonts,graphicx}

\numberwithin{equation}{section}

\title{Simple Methods for Finding Actual Congruent Number Triangles}

\author{Allan J. MacLeod,\\Statistics, O.R. and Mathematics Group (retired),\\
University of the West of Scotland,\\High St., Paisley,\\Scotland.  PA1 2BE\\
(e-mail: peediejenn@hotmail.com)}

\date{}

\begin{document}

\maketitle

\begin{abstract}
We present a, hopefully, elementary mathematical treatment of the computational aspects of congruent numbers, such that
an amateur could understand the problem and perform their own calculations.
\end{abstract}

\newpage

\section{Introduction}
A {\bf congruent number} is a positive integer which is the area of a right-angled triangle with rational sides.
Consider the right-angled triangle

\begin{center}
\begin{picture}(10,8)
\put(0,0){\line(1,0){12}}
\put(12,0){\line(0,1){8}}
\put(0,0){\line(3,2){12}}
\put(12,0.5){\line(-1,0){0.5}}
\put(11.5,0){\line(0,1){0.5}}
\put(6,0.25){b}
\put(12.25,4){a}
\put(6,4.5){h}
\end{picture}
\end{center}

By Pythagoras, $a^2+b^2=h^2$ and the area $\Delta=ab/2$. The most well-known such triangle is probably the $(3,4,5)$ triangle,
introduced at school. It has area $6$. A natural question is {\bf what other integers can occur as the area of a right-angled
triangle with integer sides?} A follow-up question would be {\bf for such integers, can we give explicit values for the sides?}

If we scale the $(3,4,5)$ sides by $2$, we get $(6,8,10)$ with $6^2+8^2=10^2$ and area $24$. Scale by $3$ and the area goes up by
factor of $9$ to $54$. Thus, we can get any square multiple of $6$ by scaling the $(3,4,5)$ triangle by an integer factor.

The standard parametric form of Pythagorean triples, with no common factor, is $a=p^2-q^2, b=2p\,q, h=p^2+q^2$ with
$p,q \in \mathbb{Z}$, $\gcd(p,q)=1$ and $p,q$ of opposite parities, which allows us to set up the simple Table $1$.

\begin{table}
\begin{center}
\caption{Basic Areas}
\begin{tabular}{lrr}
$\,$&$\,$&$\,$\\
$(p,q)$&$(a,b,h)$&$N$\\
$\,$&$\,$&$\,$\\
$(2,1)$&$(3,4,5)$&$6$\\
$(3,2)$&$(5,12,13)$&$30$\\
$(4,1)$&$(15,8,17)$&$60$\\
$(4,3)$&$(7,24,25)$&$84$\\
$(5,2)$&$(21,20,29)$&$210$
\end{tabular}
\end{center}
\end{table}

Thus, we can say $N=pq(p^2-q^2)$ is the area of a right-angled triangle with integer sides, but they seem very sparse.
What if we allow $a,b,c$ to be rational? We can then scale the $(4,1)$ triangle by $1/2$ to give $(15/2,4,17/2)$ with an area of $15$.
Similarly, the $(7,24,25)$ triangle can be scaled down to $(7/2,12,25/2)$ with area $21$. This scaling, upwards or downwards, means that we only need
to look at those $N$ which are squarefree (no square factors), and we assume this from now on.

Let $\Delta=5$ and $a=A/D$, $b=B/D$ and $h=H/D$ with $A,B,H,D \in \mathbb{Z}$. Then
\begin{equation*}
A\,B=10D^2 \hspace{2cm} A^2+B^2=H^2
\end{equation*}
so that we require
\begin{equation*}
A^4+100D^4=A^2H^2
\end{equation*}

Searching for positive $(A,D)$ pairs which give $A^4+100D^4$ equal to a square, which is divisible by $A$, quickly gives $A=9$ and $D=6$, giving
$a=3/2, b=20/3, h=41/6$. Applying this method for squarefree $N<100$, with $A+D<10000$ finds solutions for
\begin{equation*}
N=5,6,7,14,15,21,22,30,34,39,41,46,65,70,78,85
\end{equation*}

The basic problem can also be phrased differently.

\textbf{Lemma} There is a right-angled triangle with rational sides and area $N$ if and only if there are $3$ squares forming an arithmetic
progression (AP) with common difference $N$.

\textbf{Proof}
Let $(a,b,h)$ be the sides of such a triangle. Then $(a+b)^2=h^2+4N$ and $(a-b)^2=h^2-4N$, so, if we define $p=(a-b)/2$,$q=h/2$ and
$r=(a+b)/2$, we have $p^2,q^2,r^2$ is an arithmetic progression with common difference $N$.

Conversely, let $p^2,q^2,r^2$ be such an AP. Thus $r^2=q^2+N$ and $p^2=q^2-N$, and $2N=r^2-p^2=(r+p)(r-p)$. Define $a=r+p, b=r-p$,
giving $a\,b/2=N$ and $a^2+b^2=4q^2=(2q)^2$, so define $h=2q$. $\, \Box$

\vspace{0.25cm}

It is the AP property that gives this problem its name, which is derived from the Latin word for such a progression \textbf{congruum}.
For information on the history of the problem until the early 20th century, look at Chapter XVI of Dickson \cite{dick1}.
This question has been around for over $1000$ years and is still not completely solved. There a several extremely advanced arguments which can
suggest which $N$ are congruent, but no totally conclusive theorem. Computationally, we can prove a value $N$ is congruent by
computing a specimen triangle, and this is our goal - to compute an actual triangle.

There is, thus, a dual problem.

For $\Delta=N$, does a solution exist? If it does, what is it?

Practically all work has been done on the first question.
For example, the following are congruent numbers - primes of the form $8K+5$ or $8K+7$, whilst $N$ is not
a congruent number if $N$ is a prime of the form $8K+3$ or of the form $2p$ with $p$ prime and $p=16K+9$. There are a plethora of other
conditions, usually involving primes, see Chandrasekar \cite{chand} for some further examples.

We concentrate on the second problem, which is certainly not as trivial as some think. Very advanced methods can be used, such as Elkies'
wonderful Heegner-point method \cite{elk1}, but we will concentrate on as simple methods as possible.

If you are interested in pursuing this subject, you will need
a computer and some software. Any modern computer should do, and the software should be able to handle very large integers. I use the
package Pari-gp \cite{PARI2} which is excellent, easy-to-use and free!

\section{Hartley's method}
This method is over $200$ years old but provides a simple introduction to some of the ideas which appear in the rest of this report.
It is described on page $464$ of Dickson for the specific case of $N=13$. The original appeared in $1803$ in a supplement to the Ladies Diary
magazine printed in London. It is easy to generalize to $N$, but, as we will see, $N$ is restricted to a subset of possible values.

The method involves looking for an AP of squares $(x^2-N,x^2,x^2+N)$. Set $x^2+N=(x+y)^2$ and $x^2-N=(x-y\,z)^2$. Then simple algebra gives
\begin{equation*}
x=\frac{N-y^2}{2y}=\frac{N+y^2z^2}{2yz}
\end{equation*}
so that $y^2=N(z-1)/(z(z+1))$.

Hartley then looks for a value of $z$ of the form $(r^2+s^2)/(2rs)$ with $r,s \in \mathbb{Z}$. This would give
\begin{equation*}
y^2=\frac{2Nrs(r-s)^2}{(r^2+s^2)(r+s)^2}
\end{equation*}

To satisfy this equation, Hartley sets $r^2+s^2=N$ and $2rs=\Box$. This is done without any background theory, just trying to see if it works.
The first relation, thus, restricts possible $N$ to those integers which can be
written as the sum of $2$ squares, such as $5,13,29$. Let $r_0,s_0$ be such a representation.

Hartley defines $r=r_0-g\,t$ and $s=s_0-t$, where $g,t$ are rational. Thus, $t=2(g\,r_0+s_0)/(g^2+1)$, and substituting into $2\,r\,s=\Box$, gives
that $g$ must satisfy the quartic
\begin{equation}\label{hartq4}
w^2=-2(r_0g^2+2s_0g-r_0)(s_0g^2-2r_0g-s_0)
\end{equation}

The requirement that a quartic be made square is another common feature of most methods.
It is easy to write a simple search for a quartic to be a square. Unfortunately, there are two major problems.
Firstly, the quartic might not give {\bf any} square values, and secondly, it might take a long time to find one, even if it exists.

The test for whether a solution exists is quite advanced. It is possible to avoid it at the risk of considering quartics with no solutions. Searching
can be done very efficiently, up to a reasonable level, with the Pari subroutine {\bf hyperellratpoints} which is an implementation of Michael Stoll's
excellent ratpoints code.

As an example, consider $N=37$ with $r_0=6, s_0=1$. The quartic is
\begin{equation*}
w^2=-4(g^2-12g-1)(3g^2+g-3)
\end{equation*}
and it is straightforward evaluation to show that $(12, \pm 42)$ are solutions.

This gives $t=146/145$, $r=-882/145$ and $s=-1/145$. Thus $z=777925/1764$ and $y=37002/128035$. Finally $x=605170417321/9475102140$. All this,
using Lemma $1$, gives a triangle with $a=777923/6090$ and $b=450660/777923$.

The method can also be applied to $N=101=10^2+1^2$ and $N=157=11^2+6^2$. $N=101$ was only solved by Bastien \cite{bast} in $1915$, whilst $N=157$ was solved
by Don Zagier in the $1970s$. We find $g=149/141$ gives a solution for $N=101$,
but $g=-186067/136146$ is needed for $N=157$. This would take a long time for a simple search.

We can speed up such a simple search as follows.
The first thing to note is that if $g=u/v$ gives a solution, then $g=-v/u$ also gives a solution, so we need only search over $|g| \le 1$.
Secondly, the quartic has $4$ real roots, namely
\begin{equation*}
\frac{-s_0 \pm \sqrt{N}}{r_0} \hspace{2cm} \frac{r_0 \pm \sqrt{N}}{s_0}
\end{equation*}
so we have simple bounds for the quartic to be positive.

We can perform a further operation, known as a descent.
The quartic \eqref{hartq4} is of the form
\begin{equation}
(a\,x^2+b\,x+c)(d\,x^2+e\,x+f)=w^2
\end{equation}
which we consider as
\begin{equation}
a\,x^2+b\,x+c=k\,y^2  \hspace{2cm} d\,x^2+e\,x+f=k\,z^2
\end{equation}
where $k$ is squarefree.

Suppose $(x_0,y_0)$ is a rational solution of the first quadric (as such functions are called). Then, the line $y=y_0+t(x-x_0)$ will meet the quadric
at one further point, with
\begin{equation}
x=\frac{a\,x_0+b+k\,t(t\,x_0-2y_0)}{k\,t^2-a}
\end{equation}

Substitute this into the second quadric. Let $x_0=p_0/q_0$ and $y_0=w_0/q_0$ with $p_0,q_0,w_0 \in \mathbb{Z}$. Then the following
integer coefficient quartic, in $t$, must be made square, if possible.
\begin{equation*}
k^3(dp_0^2+ep_0q_0+f\,q_0^2)t^4-2k^3w_0(2d\,p_0+e\,q_0)t^3+
\end{equation*}
\begin{equation}
k^2(2a(d\,p_0^2-f\,q_0^2)+b\,q_0(2d\,p_0+e\,q_0)+4d\,k\,w_0^2)t^2-
\end{equation}
\begin{equation*}
2k^2w_0(a(2d\,p_0-e\,q_0)+2b\,d\,q_0)t+
\end{equation*}
\begin{equation*}
k(a^2(d\,p_0^2-e\,p_0q_0+f\,q_0^2)+a\,b\,q_0(2d\,p_0-e\,q_0)+b^2d\,q_0^2)
\end{equation*}

Of course, we have not yet said what values to use for $k$. These come from the resultant of the two quadratics, which is the determinant
of the matrix
\begin{equation}
\left( \begin{array}{rrrr}r_0&2s_0&-r_0&0\\0&r_0&2s_0&-r_0\\-2s_0&4r_0&2s_0&0\\0&-2s_0&4r_0&2s_0 \end{array} \right)
\end{equation}
which can be computed to be $-16(r_0^2+s_0^2)^2=-16N^2$. Thus $k$ is a squarefree divisor of $2N$, positive or negative.

This extra descent finds a solution of the $N=157$ triangle with $t=-52/131$.
We find a solution for $N=1093$ with the numerator and denominator
of both altitude and base of the triangle having roughly $50$ decimal digits.

\section{Elliptic Curve Formulation}
We now show the link between the triangle problem and elliptic curves. There are many approaches, and I particularly like this one, from
Keith Conrad \cite{con}.

Let $a,b,h > 0$ with $a^2+b^2=h^2$ and $ab=2N$, with $N \in \mathbb{Z}$.

Then $t=h-a>0$, and
\begin{equation*}
2\,a\,t=b^2-t^2 \hspace{2cm} \frac{4Nt}{b}=b^2-t^2
\end{equation*}

Multiply by $N^3b/t^3$ giving
\begin{equation*}
\frac{4N^4}{t^2}=\frac{N^3b^3}{t^3}-N^2\frac{Nb}{t}
\end{equation*}

Now, define $Y=2N^2/t$ and $X=Nb/t$ so that
\begin{equation}\label{econg}
Y^2=X^3-N^2X=X(X-N)(X+N)
\end{equation}
and we have $3$ finite rational points on the curve, namely $(0,0)$, $(N,0)$, and $(-N,0)$.

Conversely, let $(u,\pm v)$ be rational points on the elliptic curve with $|u\,v|>0$. Define $a=v/u$ choosing the sign
of $v$ to make $a>0$. Define $b=2N/a$, so $a,b$ are strictly positive with $N=a\,b/2$ and
\begin{equation*}
a^2+b^2=\frac{v^2}{u^2}+\frac{4N^2u^2}{v^2}=\frac{v^4+4N^2u^4}{u^2v^2}
\end{equation*}

Using $v^2=u^3-N^2u$, we have
\begin{equation*}
a^2+b^2=\frac{(u^3-Nu)^2+4N^2u^4}{u^2v^2}=\frac{(u^3+N^2u)^2}{u^2v^2}=\left( \frac{u^2+N^2}{v} \right)^2
\end{equation*}
so we can set $h=|(u^2+N^2)/v|$.

\includegraphics{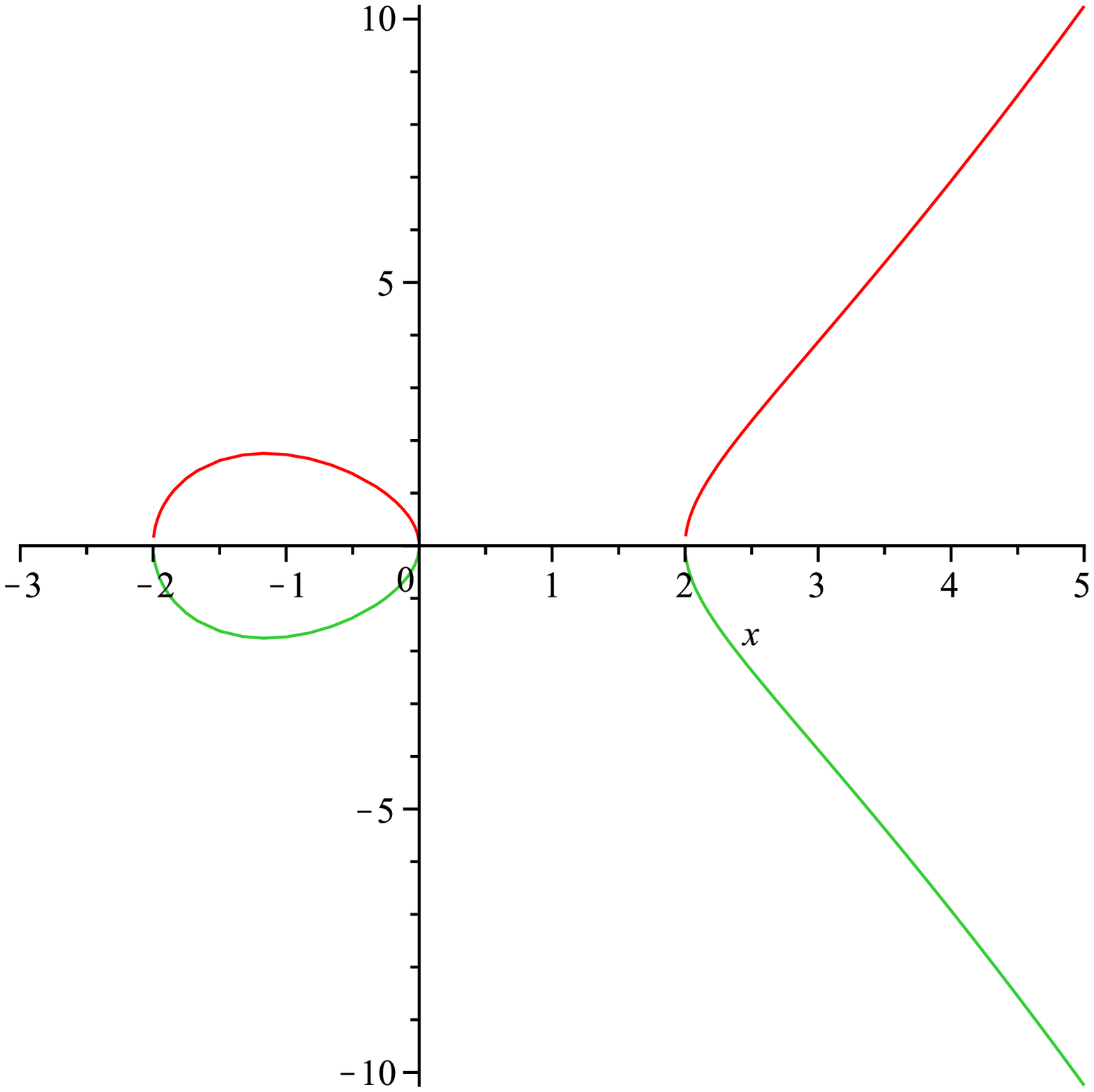}

As an example, the graph of $y^2=x^3-4x$ is shown. All congruent number elliptic curves have this basic shape.
There are $2$ components, a closed convex egg-shape for $-N \le x \le 0$ and an infinite component for $x\ge N$.

Other transformations are possible. We have
\begin{equation*}
(a+b)^2=h^2+4N \hspace{2cm} (a-b)^2=h^2-4N
\end{equation*}
and thus
\begin{equation*}
\left(\frac{h^2}{4} - N \right) \frac{h^2}{4} \left( \frac{h^2}{4} + N \right)=
\frac{h^2}{4} \frac{(a^2-b^2)^2}{16}=\Box
\end{equation*}
so $h^2/4$ gives an x-coordinate of a point on $Y^2=X^3-N^2X$.

The problem is that this x-coordinate is a square, whereas many curves have rational points which are non-squares. For example,
the $(3,4,5)$ triangle gives the point $(12, \pm 36)$.

\section{Geometry of the elliptic curve}
For $N=6$, we know some points on $Y^2=X^3-36X$, namely $(0,0)$, $(-6,0)$, $(6,0)$. The $(3,4,5)$ triangle gives $(12, \pm 36)$.
We can generate new points by simple geometrical calculations.

Suppose we join $(0,0)$ by a straight line to $(12,36)$. The line has gradient $3$ and can be expressed as $Y=36+3(X-12)$, and will meet the curve where
\begin{equation*}
X^3-36X-(36+3(X-12))^2=0 \Rightarrow X^3-9X^2-36X=0=X(X-12)(X+3)
\end{equation*}
so the line meets the curve also at $X=-3$ where $Y=-9$. By symmetry, $(-3,9)$ also lies on the curve. These points both give the $(3,4,5)$ triangle.

Joining $(12,36)$ to $(-6,0)$ gives the points $(-2, \pm 8)$, whilst joining to $(6,0)$ gives $(18,\pm 72)$, and these lead to a $(4,3,5)$ triangle.

We can also look at the tangent to the curve at a point. We have
\begin{equation*}
\frac{dY}{dX}=\frac{3X^2-36}{2Y}
\end{equation*}
so, at $(12,36)$ the tangent will have equation $Y=36+11(X-12)/2$ and meets the curve where $(X-12)^2(4X-25)=0$. The double root comes from being
a tangent, so there is only one further point of intersection at $(25/4, 35/8)$. This point gives a right-angled triangle with $a=7/10$, $b=120/7$
and $h=1201/70$ but still with area $6$.

These chord and tangent methods can be continued forever to produce an infinite sequence of right-angled triangles all of area $6$, but where
the numbers involved get larger and larger.

\section{2-isogeny Curve}
There is another related elliptic curve in this problem. Suppose $a^2+b^2=h^2$ and $ab/2=N$, then
\begin{equation*}
a^2+\frac{4N^2}{a^2}=h^2  \Rightarrow a^6+4N^2a^2=h^2a^4
\end{equation*}
so if we define $U=a^2$ and $V=ha^2$ we have the elliptic curve
\begin{equation}\label{isog}
V^2=U^3+4N^2U
\end{equation}

Now, let $(P,Q)$ be a point on \eqref{econg} with $P \ne 0$. Then, if
\begin{equation}
U=\frac{Q^2}{P^2} \hspace{2cm} V=\frac{Q(P^2+N^2)}{P^2}
\end{equation}
we have a rational point on \eqref{isog}.

In fact, we can go in the other direction. If $(U,V), U \ne 0$ is a rational point on \eqref{isog}, defining
\begin{equation}
X=\frac{V^2}{4U^2} \hspace{2cm} Y=\frac{V(U^2-4N^2)}{8U^2}
\end{equation}
gives a rational point on \eqref{econg}.
The transformations show that $X$ will have roughly double the number of decimal digits as $U$.

For example, if $N=53$, it is easy to find a point $U=4335188/14161$ $V=9552116860/1685159$ on $V^2=U^3+4*53^2U$.
This gives
\begin{equation*}
X=99278557225/1158313156 \hspace{0.5cm} Y=24583549770420915/39422029951304
\end{equation*}
 on the congruent number curve $Y^2=X^3-53^2X$. The increase in the size of the numbers involved is clear.

 It is, occasionally, possible to find a point with smaller size numbers by using the tangent to a point method discussed in the previous section.

 Let $(P,Q)$ be a rational point on \eqref{econg}. The tangent to this point has gradient $(3P^2-N^2)/(2Q)$ so the line
 \begin{equation*}
 Y=Q+\frac{3P^2-N^2}{2Q} \, (X-P)
 \end{equation*}
meets \eqref{econg} at only one further point, which standard algebra gives as
\begin{equation}
X=\frac{(P^2+N^2)^2}{4P(P+N)(P-N)} = \frac{(P^2+N^2)^2}{4Q^2}
\end{equation}
which is a square, and, sometimes, we can solve the quartic
\begin{equation}\label{tanq4}
(P^2+N^2)^2-4XP(P-N)(P+N)
\end{equation}

For example, $N=62$ gives the following data. The point
\begin{equation*}
(7150393600/32798529, 695599219282240/187837175583)
\end{equation*}
lies on $V^2=U^3+4*62^2U$ so that
\begin{equation*}
X=4229297547568411201/58630597962753600
\end{equation*}
gives a rational point on $Y^2=X^3-62^2X$.

The quartic \eqref{tanq4} factors into $4$ linear terns giving possible $P$ values
\begin{equation*}
\frac{124002}{529} \hspace{1cm} \frac{-706831}{19600} \hspace{1cm} \frac{2430400}{22801} \hspace{1cm} \frac{-1016738}{62001}
\end{equation*}
all of which are small.

The two curves are said (in the parlance of elliptic curves) to be 2-isogenous to each other.

\section{Basic Methods to Find Triangles}
All the methods are based on the fact that rational points on \eqref{econg} are such that
\begin{equation*}
X=\frac{d\,u^2}{v^2} \hspace{2cm} Y=\frac{d\,u\,w}{v^3}
\end{equation*}
with $d,u,v,w \in \mathbb{Z}$, $d$ is squarefree, and $\gcd(u,v)=\gcd(w,v)=\gcd(d,v)=1$.

Putting these into \eqref{econg} gives (after cancellation)
\begin{equation}\label{uvw}
d^2w^2=d^3u^4-N^2d\,v^4
\end{equation}

Dividing by $d^2$ shows that $d|N$, so we have a (usually) small set of possible $d$ values, including $\pm 1$ and $\pm N$. If $N$
is prime these are the only possible values.

First try $d=1$, so we look for solutions of
\begin{equation*}
w^2=u^4-N^2v^4
\end{equation*}
and we just search with $u,v$ positive and $\gcd(u,v)=1$, with $u+v \le L$, looking for a square.
We can also search for $d=N$ using
\begin{equation*}
w^2=N(u^4-v^4)
\end{equation*}

Using $L=9999$, we found the following extra congruent numbers up to $N=97$.
\begin{equation*}
13,23,29,31,37,47,55,61,71,87,95
\end{equation*}
though we often get the triangles with $a,b$ reversed and, for $N=41$, we get two completely different triangles.

If we look closely at the results (always a good idea), we find the first $12$ congruent numbers found are
\begin{equation*}
5,6,7,13,14,15,21,22,23,29,30,31
\end{equation*}
and these are all of the form $8K+5$, $8K+6$ or $8K+7$, and it is now a famous conjecture that this is true, though the proof is
taking some time.

It would be nice at this point to have a better idea of how many congruent numbers there are in $[0,99]$. Tunnell's results
provide a simple computational criterion, though based on the truth of the Birch and Swinnerton-Dyer (BSD) conjecture. You don't
need to understand what is going on, just be able to apply it.

For a given square-free integer $N$, define
\begin{equation*}
A_N = \# \{ ( x , y , z ) \in \mathbb{Z}^3 : N = 2 x^2 + y^2 + 32 z^2 \}
\end{equation*}
\begin{equation*}
B_N = \# \{ ( x , y , z ) \in \mathbb{Z}^3 : N = 2 x^2 + y^2 + 8 z^2 \}
\end{equation*}
\begin{equation*}
C_N = \# \{ ( x , y , z ) \in \mathbb{Z}^3 : N = 8 x^2 + 2 y^2 + 64 z^2 \}
\end{equation*}
\begin{equation*}
 D_N = \# \{ ( x , y , z ) \in \mathbb{Z}^3 : N = 8 x^2 + 2 y^2 + 16 z^2 \}
\end{equation*}
where $\#S$ denotes the number of elements in the set $S$. You just have to be careful about what happens if one of $x,y,z$ is zero.

Tunnell's theorem states that supposing $N$ is a congruent number, $2A_N=B_N$ if $N$ is odd, and $2C_N=D_N$ if $N$ is even. In the opposite
direction, if the BSD Conjecture is true these set-size equalities are enough to prove $N$ is congruent or not.

We thus find that we still need the following values
\begin{equation*}
38, 53, 62, 69, 77, 79, 86, 93, 94
\end{equation*}

Hartley's method disposes of $N=53$ giving the triangle
\begin{equation*}
a=\frac{1472112483}{202332130} \hspace{2cm} b=\frac{21447205780}{1472112483}
\end{equation*}

All but $53$ and $79$ in these values are composite, so we should include information on possible factors. In \eqref{uvw}, set
$N=d*e$ giving
\begin{equation*}
w^2=d(u^4-e^2v^4)
\end{equation*}
so that $d|w$. Set $w=dz$, so we look for solutions of $d z^2=u^4-e^2v^4$, just by looping round factors of $N$. This finds
solutions for these values quickly, except for $N=79$.

\section{Descent Methods}
In this section, we give a further example of a descent - trying to reduce the problem to one with
smaller numbers in the solution. The congruent number problem has a large number of possible descents, with which the reader
is encouraged to experiment.

Given the basic identity,
\begin{equation*}
d^2w^2=d^3u^4-N^2dv^4
\end{equation*}
we use the fact that $d|N$ to write $N=de$ with $e$ also an integer, leading to
\begin{equation*}
w^2=d(u^4-e^2v^4)
\end{equation*}

Thus, $d|w$ so we can write $w=dz$ with
\begin{equation}
dz^2=u^4-e^2v^4=(u+ev^2)(u-ev^2)
\end{equation}

Set $u^2+ev^2=dp^2$ and $u^2-ev^2=q^2$, so $z=p\,q$. Thus $ev^2=u^2-q^2=(u+q)(u-q)$, so set $u+q=e\,r^2$, $u-q=s^2$. We have
$v=r\,s$ and $u=(er^2+s^2)/2$, with
\begin{equation*}
\frac{(er^2+s^2)^2}{4}+er^2s^2=d \, \Box
\end{equation*}
and finally the quartic
\begin{equation}\label{eqdesc}
e^2r^4+6e\,r^2s^2+s^4=4\,d\,t^2
\end{equation}

As we mentioned in the previous section, we still need a solution for $N=79$. $N=79$ uses $d=-79,e=-1$ needing $r=125, s=52$ for the solution,
$a=233126551/167973000$ and $b=335946000/2950969$.

As a further example, consider $N=101$, which has $d=\pm 1$ or $d=\pm N$. Trying $d=-1$ so $e=-101$ gives the quartic
\begin{equation*}
10201r^4-606r^2s^2+s^4=- \, \Box
\end{equation*}
and it is a simple search to find $r=53, s=397$. These give $v=21041$ and $u=-63050$ leading to
\begin{equation*}
a=\frac{44538033219}{1326635050} \hspace{0.75cm} b=\frac{267980280100}{44538033219} \hspace{0.75cm} h=\frac{2015242462949760001961}{59085715926389725950}
\end{equation*}

The astute reader will have noticed that none of these substitutions have been proven. If they work - they work !

\section{A Special Descent}
The result that is almost the "poster-child" for congruent numbers is Don Zagier's solution for $N=157$, which is always given to show the sizes
of the numbers involved. Nobody describes how it was computed and, in fact, some get it wrong, ascribing the solution to the use of Heegner points - a
subject on which Zagier was one of the first experts, see Dalawat \cite{dal}. In a reply, on 10 Nov. 2015, to a question on the web-site MathOverflow,
Noam Elkies stated that it was actually a descent argument which Zagier had used. Carlo Beenakker gave the reference \cite{zag} and a translation from
the original German to English.

This description was extremely sparse, so I set out to try to work out a possible simple descent, which I now present. In \eqref{uvw}, set $d=1$
giving
\begin{equation}
w^2=u^4-N^2v^4=(u^2+Nv^2)(u^2-Nv^2)
\end{equation}

Now set $u^2+Nv^2=r\,s^2$ and $u^2-Nv^2=r\,t^2$, with $r$ squarefree. Thus, $2u^2=r(s^2+t^2)$, $2Nv^2=r(s^2-t^2)$ and $w=s\,t$.
Assume $N$ is prime and let $r=N$. So, $N|(2u^2)$, so, since $N>2$, we have $N|u$. Set $u=Nz$, with $2Nz^2=(s^2+t^2)$.

Also, $2v^2=s^2-t^2$ which we write $t^2+2v^2=s^2$.
We can derive a parametric solution to this last identity as follows. Let $X=v/t, Y=s/t$ so $Y^2=1+2X^2$, which has a trivial solution
$X=0, Y=1$. The line $Y=1+kX$ passes through this point and meets this curve again at one point. This gives $t=2p^2-q^2$, $v=2pq$ and $s=2p^2+q^2$.

We end up with $Nz^2=4p^4+q^4=(2p^2)^2+(q^2)^2$.

We now assume $N$ is prime of the form $8K+5$, so $N$ can be expressed as the sum of two integer squares, eg. $157=11^2+6^2=6^2+11^2$. Thus $N=c^2+d^2$,
and assume $z=x^2+y^2$. Thus
\begin{equation}
(c^2+d^2)(x^2+y^2)^2=(2p^2)^2+(q^2)^2
\end{equation}
which can be written, in terms of complex numbers
\begin{equation}
(c+di)(c-di)(x+yi)^2(x-yi)^2=((2p^2)+(q^2)i)((2p^2)-(q^2)i)
\end{equation}

We now have
\begin{equation}
(c+di)(x+yi)^2=2p^2+q^2i
\end{equation}
which gives, by equating the real and imaginary parts
\begin{equation}
c(x^2-y^2)-2dxy=2p^2 \hspace{2cm} d(x^2-y^2)+2cxy=q^2
\end{equation}

If we find a simple solution $(x_0,y_0,p_0)$ to the first equation, we can parameterize, as before, giving
\begin{equation*}
\frac{x}{y}=\frac{cx_0-2(dy_0-k(kx_0-2p_0))}{y_0(2k^2-c)}
\end{equation*}
where $k$ is a rational parameter.

Substituting into the second equations, and simplifying the result, gives the quartic
\begin{equation*}
j^2=4(2cx_0y_0+d(x_0^2-y_0^2))k^4-16p_0(cy_0+dx_0)k^3+
\end{equation*}
\begin{equation}
4d(c(x_0^2-y_0^2)-2dx_0y_0+4p_0^2)k^2+8p_0(c^2y_0-cdx_0+2d^2y_0)k-
\end{equation}
\begin{equation*}
2c^3x_0y_0+c^2d(x_0^2+3y_0^2)-4cd^2x_0y_0+4d^3y_0^2
\end{equation*}

Now, consider $N=157$. There is no solution $(x_0,y_0,p_0)$ for $c=11,d=6$, but $(4,1,1)$ is a solution when $c=6,d=11$.
This gives the quartic
\begin{equation*}
j^2=4(213k^4-200k^3+66k^2+28k-124)
\end{equation*}
and a simple search finds a solution when $k=-262/79$. This gives $x/y=322213/49921$, $p=356441$ and $q=1143522$.

From this, we have $v=815196250404$, together with $s=1561742937446$ and $t=-1053542191522$.
Finally $u=16691191806770$, and we have the point
\begin{equation*}
\left ( \, \frac{16691191806770^2}{815196250404^2} \, \, , \, \, \frac{538962435089604615078004307258785218335}{815196250404^3} \, \right )
\end{equation*}
on the curve $Y^2=X^3-157^2X$.

As a measure of the size of the numbers found, we use the {\bf height} of a rational point on an elliptic curve. Very, very roughly,
the height gives an indication of the number of decimal digits in the numerator and denominator of the X-coordinate of a point.
For example, the rational point for the $N=157$ congruent number elliptic curve has height $54.6$, with the numerator having $24$ digits and the denominator $22$.

Experiments find the largest height point found so far, with this method, to be for $N=7309$ with height $121.05$.

\section{A 2-isogeny descent}
In the section, we describe a descent on the 2-isogenous curve $Y^2=X^3+4N^2X$. It is based on the observation that the X-coordinate
of a rational point can often be of the form $2u^2/v^2$. Thus, we suppose $X=2u^2/v^2$ and $Y=2uw/v^3$, so that
\begin{equation*}
w^2=2(u^4+n^2v^4)
\end{equation*}
and so $2|w$, giving $w=2t$ and
\begin{equation*}
2t^2=u^4+N^2v^4
\end{equation*}

This shows much fewer factorisation possibilities than previous sections. By employing complex numbers, we can proceed. First note
that $2=(1+i)(1-i)$ and assume $t=a^2+b^2=(a+bi)(a-bi)$. Thus
\begin{equation*}
(1+i)(1-i)(a+bi)^2(a-bi)^2=(u^2+Nv^2i)(u^2-Nv^2i)
\end{equation*}

Define
\begin{equation*}
u^2+Nv^2i=(1+i)(a+bi)^2 \hspace{2cm} u^2-Nv^2i=(1-i)(a-bi)^2
\end{equation*}
giving
\begin{equation}
u^2=a^2-2ab-b^2 \hspace{2cm} Nv^2=a^2+2ab-b^2
\end{equation}

The first quadric $u^2=a^2-2ab-b^2$ can be parameterized in the standard way to give
\begin{equation*}
\frac{b}{a}=\frac{-2(k+1)}{k^2+1}
\end{equation*}
with $k \in \mathbb{Q}$. This means that
\begin{equation*}
\frac{k^4-4k^3-6k^2-12k-7}{N}=\Box
\end{equation*}

As an example, consider $N=103$. Note that the algebra is exactly the same if $N$ is negative, and we find $k=100/19$ gives a solution to the last equation for $N=-103$. This gives $b/a=-4522/10361$.

Using the numerator and denominator, $u=13439$ and $v=257$ and the point
\begin{equation*}
\left ( \frac{361213442}{66049} , \frac{6869952561580}{16974593} \right )
\end{equation*}
on $V^2=U^3+4*103^2U$.

This then gives the point
\begin{equation*}
\left( \frac{16332534559428025}{11928893315329}, \frac{2081363685506152106939880}{41200286097029552767} \right)
\end{equation*}
on $Y^2=X^3-103^2X$.

\section{Solution of \eqref{eqdesc}}
We now return to solving equation \eqref{eqdesc}, based on the paper by Komoto, Watanabe and Wada \cite{kww}.
I generalized this in a recent arXiv preprint \cite{mac2}.

Define $r^2=f$ and $s^2=g$, so the basic equation is
\begin{equation*}
e^2f^2+6e\,f\,g+g^2=4\,d\, t^2
\end{equation*}
which we write as
\begin{equation}
(f,g,t) \left(\begin{array}{lrr}e^2&3e&0\\3e&1&0\\0&0&-4\,d\end{array} \right) \left( \begin{array}{r}f\\g\\t \end{array} \right) = 0
\end{equation}

Let $(1,g_0,t_0)$ be an initial solution of this quadratic. Setting $f=1$ makes things a lot easier later on, though we cannot always find such a solution.
Define
\begin{equation}
\left( \begin{array}{l}f\\g\\t \end{array} \right) = \left( \begin{array}{lrr}1&0&0\\g_0&1&0\\t_0&0&1 \end{array} \right) \left( \begin{array}{l}h\\j\\k \end{array} \right)
\end{equation}
so that the basic equation is now
\begin{equation}
(h,j,k) \left(\begin{array}{lrr}0&3e+g_0&-4dt_0\\3e+g_0&1&0\\-4dt_0&0&-4\,d\end{array} \right) \left( \begin{array}{r}h\\j\\k \end{array} \right) = 0
\end{equation}
or
\begin{equation*}
(6e+2g_0)hj-8dt_0hk+j^2-4dk^2=0
\end{equation*}
so that $j$ must be even. Set $j=2m$ and $3e+g_0=2\gamma$ giving
\begin{equation*}
2\gamma\,h\,m-2t_0\,d\,h\,k+m^2-d\,k^2=0
\end{equation*}
or
\begin{equation*}
(h,m,k) \left(\begin{array}{lrr}0&\gamma&-d\,t_0\\\gamma&1&0\\-d\,t_0&0&-d\end{array} \right) \left( \begin{array}{r}h\\m\\k \end{array} \right) = 0
\end{equation*}

We now suppose $\gcd(\gamma,-d\,t_0)=1$, so there exist integers $\alpha, \beta$ such that $\alpha \gamma-\beta d\,t_0=1$.
Define
\begin{equation}
\left( \begin{array}{l}h\\m\\k \end{array} \right) = \left( \begin{array}{lrr}1&0&0\\0&\alpha&d\,t_0\\0&\beta&\gamma \end{array} \right) \left( \begin{array}{l}h\\n\\p \end{array} \right)
\end{equation}
so that
\begin{equation}
2h\,n+n^2(\alpha^2-\beta^2\,d)+n\,p(2t_0\alpha\,d-2\beta\,d\,\gamma)=d\,p^2(\gamma^2-t_0^2\,d)=2\,d\,e^2\,p^2
\end{equation}

Define
\begin{equation*}
q=2h+n(\alpha^2-\beta^2\,d)+p(2t_0\alpha\,d-2\beta\,d\,\gamma)
\end{equation*}
so that
\begin{equation}
n\,q=2\,d\,e^2\,p^2
\end{equation}

We now perform a descent on this equation by setting
\begin{equation}
n=d\,e^2\,a^2 \hspace{2cm} q=2\,b^2
\end{equation}
so that $p=a\,b$.

Thus
\begin{equation*}
h=b^2+d(\beta\,\gamma-\alpha\,t_0)a\,b+d\,e^2(bet^2\,d-alp^2)a^2/2
\end{equation*}
which we write as
\begin{equation}
r^2=f=h=\left (b+\frac{d(\beta\,\gamma-\alpha\,t_0)a}{2}\right )^2+Ka^2=c^2+K\,a^2
\end{equation}
with
\begin{equation*}
K=\frac{-d(\alpha^2(d\,t_0^2+2e^2)-2d\alpha\beta\gamma\,t_0+d\beta^2(\gamma^2-2e^2))}{4}
\end{equation*}
Somewhat surprisingly, we can show $K=-d/4$.

From $r^2-c^2=K\,a^2$, we set $r+c=K\,x^2, (r-c)=y^2$, so $a=x\,y$. We also have
\begin{equation*}
s^2=g_0h+2(\alpha\,n+d\,t_0\,p)
\end{equation*}
which, after lots of substitutions, means that the following quartic must be a square
\begin{equation}
d^2g_0x^4-16d^2t_0x^3y+
\end{equation}
\begin{equation*}
8d(8\alpha(d\,t_0^2+2e^2)-8\beta d\gamma\,t_0-g_0)x^2y^2-64d\,t_0xy^3+16g_0y^4
\end{equation*}

It is straightforward to write a program to perform these calculations, and
we can easily experiment. We find no solutions if $N=8M+1,8M+7$. For other $N$, we find heights of up to about $110$, if we use hyperellratpoints with
a search limit of $99999$. For $N=8M+6$, we find that some of the quartics factorize further into $2$ quadratics, and we can then use the descent
described in section $2$.

For example, $N=8662$ gives a quartic (setting $x/y=z$)
\begin{equation*}
\Box=405z^4-25444z^3+217608z^2-101776z+6480
\end{equation*}
with the quartic factorising to
\begin{equation*}
\Box=(5z^2-264z+20)(81z^2-812z+324)
\end{equation*}
which has a solution, from the extra descent
\begin{equation*}
z=\frac{420742975301}{561839195}
\end{equation*}
which would be very time-consuming to find from the original quartic. The corresponding point on the elliptic curve has height $231.4$, which
is the largest height found so far.

\section{Miscellaneous Calculations}
In \cite{wiman}, Wiman found that $N=1254=2*3*11*19$ gave a curve of rank $5$ and $N=29274=2*3*7*17*41$ gave rank $6$. Readers should note
that Wiman defines the rank as including the number of generators of the torsion subgroup, so is $2$ more than the modern definition. If we
return to the Pythagorean method of the Introduction, we must solve
\begin{equation*}
p\,q(p^2-q^2)=N\,D^2
\end{equation*}
and it is very easy to program this.

For $p+q \le 999$, we find $9$ different triangles for $1254$ and $10$ for $29274$. Thus, this simple method provides a good indicator
for larger rank elliptic curves.

Recently, after my main computer died, I had to ask Randall Rathbun \cite{rath} for a replacement copy of his computer files giving
solutions for $N \in [1,999999]$ (apart from $8$ rank two curves where one generator only is known). The obvious next step is to go beyond one million.

Using Tunnell's criterion, the first five values are $N=10^6+1,5,6,7,9$. The first and last values are of the form $8K+1$, and standard
conjectures would imply that there at least $2$ generators for both values. For $N=10^6+1=101*9901$, we found different solutions using
the Zagier method and Hartley method described before. These have heights $93.0$ and $101.8$ respectively.

I, however, have a personal copy of Magma, which is an extremely powerful (but not free) package. This finds a point with height $21.6$, with
x-coordinate on the elliptic curve
\begin{equation*}
\frac{310941179641}{77841}=\frac{101^2*5521^2}{3^4*31^2}
\end{equation*}
so we have $101|N$ and $101|u$, which led to the following descent.

Putting $d=1$ into \eqref{uvw} gives $w^2=u^4-N^2v^4$. Now suppose $N=a\, b$, and that $u=a\,t$, so
\begin{equation*}
w^2=a^4t^4-a^2b^2v^4=a^2(a^2t^4-b^2v^4)=a^2(at^2+bv^2)(at^2-bv^2)
\end{equation*}

Thus, we can apply the method, described at the end of the section on Hartley's method, to the quadratics
\begin{equation*}
ax^2+b=\Box \hspace{2cm} ax^2-b=\Box
\end{equation*}
which does find the solution for $N=1000001$ with height $21.6$.

As a final show of the power of these simple methods, I used Tunnell's method to give the $361$ squarefree congruent numbers in $[1,999]$ and
applied all the methods to these values. The number unsolved was reduced to $65$ in an afternoon on an old laptop. I have a general descent
algorithm for curves $y^2=x^3+Ax^2+Bx$ which is described in \cite{mac1}, which reduced the number of unsolved $N$ to $19$, one of which was
then solved by Hartley's method with a large search region. The $18$ remaining values have $15$ primes and $3$ values of twice a prime.

For the primes of the form $8M+7$, my Pari implementation of Elkies method can find a solution quickly, though this method cannot be described
as simple, so we have $343$ values where an actual triangle can be found by essentially simple methods.

\end{document}